\documentclass[12pt]{article}
\usepackage{graphicx} 

\usepackage{epsfig}
\usepackage{amsmath}
\usepackage{amssymb}
\usepackage{amsthm}
\usepackage{latexsym}

\newcommand{\comment}[1]{}
\newtheorem{theorem}{Theorem}        
\newtheorem{lemma}{Lemma}[section]

\newtheorem{corollary}{Corollary}[section]
\newtheorem{proposition}{Proposition}[section]

\newtheorem{definition}{Definition}[section]
\begin{document}
 
\title{{\LARGE\sf
A Particular Bit of Universality:
Scaling Limits of Some Dependent Percolation Models}}

\author{
{\bf Federico Camia}
\thanks{Present address: EURANDOM, P.O. Box 513, 5600 MB Eindhoven, The Netherlands.
E-mail: camia@eurandom.tue.nl}\,
\thanks{Research partially supported by the
U.S. NSF under grants DMS-98-02310 and DMS-01-02587.}\\
{\small \tt federico.camia\,@\,physics.nyu.edu}\\
{\small \sl Department 
of Physics, New York University, New York, NY 10003, USA}\\
\and
{\bf Charles M.~Newman}
\thanks{Research partially supported by the
U.S. NSF under grants DMS-98-03267 and DMS-01-04278.}\\
{\small \tt newman\,@\,courant.nyu.edu}\\
{\small \sl Courant Inst.~of Mathematical Sciences, 
New York University, New York, NY 10012, USA}\\
\and
{\bf Vladas Sidoravicius}
\thanks{Research partially supported by FAPERJ grant
E-26/151.905/2000 and CNPq.}\\
{\small \tt vladas\,@\,impa.br}\\
{\small \sl Instituto de Matematica Pura e Aplicada,
Rio de Janeiro, RJ, Brazil}\\
}

\date{}

\maketitle

\begin{abstract}
We study families of dependent site percolation models on
the triangular lattice ${\mathbb T}$ and hexagonal lattice
${\mathbb H}$ that arise by applying certain cellular automata
to independent percolation configurations.
We analyze the scaling limit of such models and
show that the distance between macroscopic portions of cluster
boundaries of any two percolation models within one of our
families goes to zero almost surely in the scaling limit.
It follows that each of these cellular automaton generated
dependent percolation models has the same scaling limit
(in the sense of Aizenman-Burchard \cite{ab})
as independent site percolation on ${\mathbb T}$.
\end{abstract}




\section{Introduction}

The phase diagrams of many physical systems have a ``critical
region'' where all traditional approximation methods, such as
mean-field theory and its generalizations, fail completely
to provide an accurate description of the system's behavior.
This is due to the existence of a ``critical point,''
approaching which, some statistical-mechanical quantities diverge,
while others stay finite but have divergent derivatives.
Experimentally, it is found that these quantities usually
behave in the critical region as a power law.
The exponents that appear in those power laws, called
{\bf critical exponents}, describe the nature of the
singularities at the critical point.
The theory of critical phenomena based on the renormalization
group suggests that statistical-mechanical systems fall into
``universality classes'' such that systems belonging to the
same universality class have the same critical exponents.

It should be emphasized that the phenomenon of universality
was discovered at least twice prior to the introduction of the
renormalization group.
In fact, in what Alan Sokal~\cite{sokal} calls the Dark Ages of
the theory of critical phenomena, before the 1940's, physicists
generally believed that all systems have the same critical exponents,
namely those of the mean-field theory of Weiss~\cite{weiss}~(1907)
or its analogue for fluids, the van der Waals theory~\cite{vanderwaals}~(1873).
The failure of this type of universality became apparent as
early as 1900, following experiments on fluid systems, but
it began to be taken seriously only after Onsager's exact
solution~\cite{onsager} (explicitly displaying non-mean-field exponents)
in 1944 and the rediscovery of the experimental evidence of non-mean-field
values for critical exponents by Guggenheim in 1945~\cite{guggenheim}.
Nonetheless, although mean-field theory is clearly incorrect
for short range models, some universality does seem to hold.
Many, if not all, different fluids, for instance, seem to have
the same value for the critical exponent $\beta$ (related to
density fluctuations), and it is believed that, for example,
carbon dioxide, xenon and the three dimensional Ising model
should all have the same critical exponents.
Maybe even more surprisingly, it was soon realized that some
(but not all) magnetic systems have the same critical exponents
as do the fluids.
This remarkable phenomenon seems to suggest the existence of
a mechanism that makes the details of the interaction irrelevant
in the critical region.
Nevertheless, the critical exponents should depend on
the dimensionality of the system and on any symmetries in the
Hamiltonian.
Despite being a very plausible and appealing heuristic idea,
backed up by renormalization-group arguments and
empirical evidence, only very few cases
are known in which universality has actually been \emph{proved}.

Percolation, with its simplicity and important physical
applications, is a natural candidate for studying
universality. This is especially so after the ground-breaking work
of Lawler, Schramm, and Werner~\cite{lsw1, lsw2, lsw3, lsw4, lsw5,
lsw6, lsw7, schramm}, who identified the only possible
conformally invariant scaling limit of critical percolation and
derived many of its properties, and that of Smirnov~\cite{smirnov, smirnov1},
who proved that, indeed, critical site percolation on the triangular
lattice has a conformally invariant scaling limit.
The combination of those results made it possible~\cite{sw} to
verify the values of the critical exponents predicted in the
physical literature in the case of critical site percolation on
the triangular lattice (and to derive also some results that had
not appeared in the physics literature, such as an analogue of
Cardy's formula ``in the bulk''~\cite{schramm1}, or the description
of the so-called backbone exponent~\cite{lsw5}).

It is generally accepted that the lattice should play no role
in the scaling limit, and that there should be no difference
between bond and site percolation.
In other words, two-dimensional critical (independent)
percolation models, both site and bond, should belong to the
same universality class,
regardless of the lattice (at least for periodic lattices like
the square, triangular or hexagonal lattice).
Once again, though, despite being a very natural and plausible
conjecture, such universality has not yet been proved.
There is however another natural direction in which to study
universality, which consists in analyzing critical percolation
models on a given lattice that differ in their dependence
structures.
It is this direction that we pursue in this paper.

The cellular automata that we use to generate our families
of dependent percolation models arise naturally in the study
of the zero-temperature limit of Glauber dynamics or as
coarsening or agreement-inducing dynamics.
The action of such cellular automata can be viewed as a sort
of ``small (local) perturbation'' of the original, independent
percolation model, and our main corollary can be viewed as
proving a form of universality for two dimensional percolation.
Therefore, we provide an explicit example of
the principle of universality, in the strong form concerning
scaling limits.

To be more precise, there are at least two, a priori
different, notions of universality, one concerning the critical
exponents discussed above, and a second one concerning the
\emph{continuum scaling limit}.
The two concepts are closely related, but in
this paper we are only concerned with the second type of
universality (the first type will be considered in a future paper).

We consider a family of dependent percolation models that arise
through a (discrete time) deterministic cellular automaton T
acting on site percolation configurations $\sigma$ on the set
of sites of the triangular lattice ${\mathbb T}$.
Each configuration $\sigma$ corresponds to an assignment
of $-1$ or $+1$ to the vertices of ${\mathbb T}$.
The variable $\sigma_x$, corresponding to the value of
$\sigma$ at $x$, is commonly called a {\bf spin} variable.
At discrete times $n=1,2,\ldots$, each spin $\sigma_x$
is updated according to the following rules
(later, in Section \ref{last}, we will introduce other cellular
automata, both on ${\mathbb T}$ and ${\mathbb H}$, generated
by different rules):
\begin{itemize}
\item if $x$ has three or more neighbors whose spin is the
same as $\sigma_x$, then the latter does not change value,
\item if $x$ has only two neighbors $y_1$ and $y_2$ that
agree with $x$, \emph{and} $y_1$ and $y_2$
are not neighbors, then $\sigma_x$ does not change value,
\item otherwise,
$\sigma_x$ changes value: $\sigma_x \rightarrow -\sigma_x$.
\end{itemize}
The starting configuration $\sigma^0$, at time $n=0$, of our cellular
automaton is chosen from a Bernoulli product measure corresponding to
independent critical percolation, and the distributions at times
$n \ge 1$ (including the final state as $n \to \infty$) of the discrete
time deterministic dynamical process $\sigma^n$ are the other members
of our family of dependent percolation models.
We show that all those dependent percolation models have the
same scaling limit, thus providing an explicit example in
which universality can be proven.
This comes as a corollary of our main result, Theorem \ref{thm}.
To explain the main result, we first need some terminology.

In the scaling limit, the \emph{microscopic scale} of the system
(i.e., the lattice spacing $\delta$) is sent to zero, while focus
is kept on features manifested on a \emph{macroscopic scale}.
In the case of percolation,
it is far from obvious how to describe such a limit and we
only do so briefly here; for more details, see \cite{a,a1,ab}.
We will make use of the approach introduced by Aizenman and
Burchard \cite{ab} (see also \cite{abnw}) applied to portions
of the boundaries between clusters of opposite sign,
and present our results in terms of closed collections of curves
in the one-point compactification $\dot{\mathbb R}^2$ of
${\mathbb R}^2$, which we identify (via the stereographic projection)
with the two-dimensional unit sphere.
For each fixed $\delta>0$, the curves are, before compactification,
polygonal paths of step size $\delta$ (i.e., polygonal paths
between sites of dual lattice).
The distance between curves is defined so that two curves are
close if they shadow each other in a metric which shrinks at
infinity (for the details, see Sections~\ref{comp}~and~\ref{curves}).

The dynamics allows one to construct the whole family of percolation
models for all $n \in \{ 0, 1, \ldots, \infty \}$ on the same probability space,
i.e., there is a natural coupling $\nu$, realized through the dynamics,
between any two percolation models in the family.
In terms of this coupling, Theorem \ref{thm} states, roughly speaking,
that the $\nu$-probability that the distance between the two collections
of curves corresponding to two distinct percolation models in the family
is bounded away from zero vanishes as $\delta \to 0$.
This means that, in the limit $\delta \to 0$, given any two models
in the family, for every curve in one of them, there exists a curve
in the other one that shadows the first curve and vice-versa.

The rest of the paper is organized as follows.
In Section~\ref{def+res}, we describe the behavior of the cellular automaton
T, define the family of dependent percolation models dynamically
generated by T, and state the main results.
The proofs of the main results are contained in Section~\ref{proofs}.
The dynamics described in Section~\ref{def+res} is chosen as a prototypical
example, but is not the only one for which our results apply.
In Section~\ref{last}, we introduce other such dynamics (both on $\mathbb T$
and $\mathbb H$), which can be obtained as suitable zero-temperature limits of
stochastic Ising models.

\section{Definitions and results} \label{def+res}

\begin{figure}[!ht]
\begin{center}
\includegraphics[width=8cm]{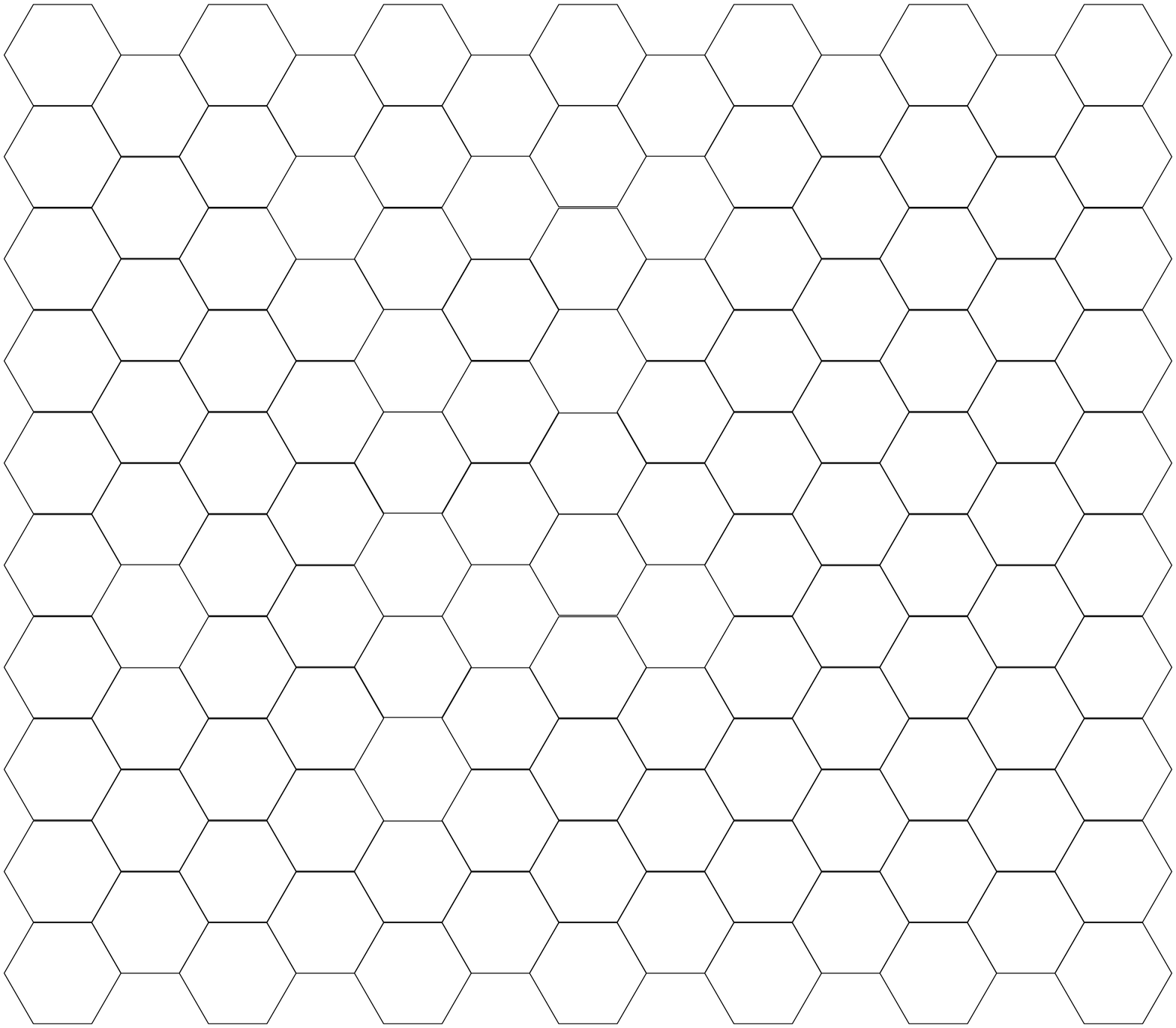}
\caption{The hexagonal (or honeycomb) lattice.}\label{Hlattice}
\end{center}
\end{figure}

We start this section by giving a more detailed definition
of one of the families of dependent percolation models
that are the object of investigation of this paper.
The models in that family are defined on the triangular lattice ${\mathbb T}$,
embedded in ${\mathbb R}^2$ by identifying the sites of ${\mathbb T}$
with the elementary cells (i.e., regular hexagons) of the hexagonal
(or honeycomb) lattice ${\mathbb H}$ (see Figure \ref{Hlattice}).
We will use those models as a paradigm and will give for them explicit
and detailed proofs of the results.
Later on, we will point out how to modify the proofs to adapt them to
the other models discussed in the paper.

In the rest of the paper,
points of ${\mathbb R}^2$ will be denoted by $u$ and $v$, while for
the sites of ${\mathbb T}$ we will use the Latin letters $x$, $y$, $z$
and the Greek letters $\zeta$ and $\xi$.
An edge of ${\mathbb T}$ incident on sites $x$ and $y$ will be denoted
by $\eta_{x,y}$, while by $\eta^*_{x,y}$ we denote the dual edge (in
${\mathbb H}$) perpendicular bisector of $\eta_{x,y}$.

\subsection{The dynamics and the percolation models}

\begin{figure}[!ht]
\begin{center}
\includegraphics[width=8cm]{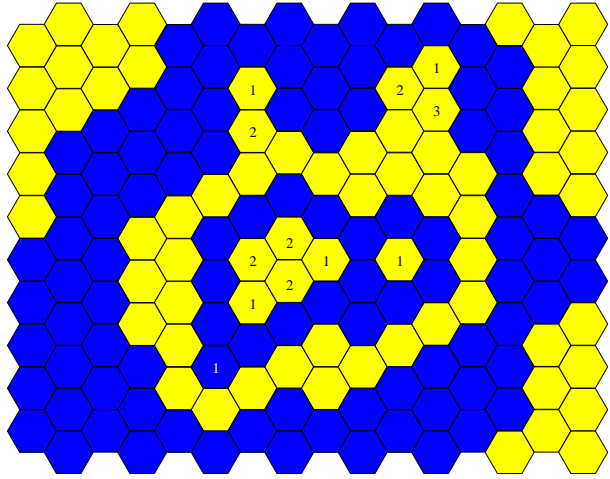}
\caption{Example of local configuration with unstable spins.
The numbered hexagons correspond to spins that will flip and the
numbers indicate at what time step the spin flips occur.
}\label{example3}
\end{center}
\end{figure}

\begin{figure}[!ht]
\begin{center}
\includegraphics[width=8cm]{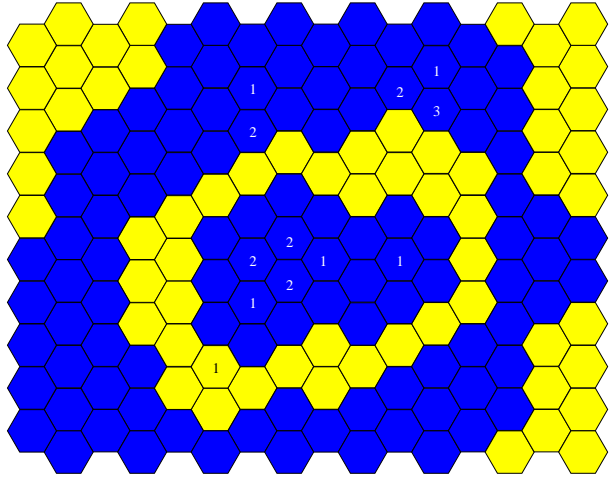}
\caption{Same local configuration after all the unstable spins have flipped.
}\label{example4}
\end{center}
\end{figure}

We construct a family of dependent percolation models by means of
a cellular automaton T acting on site percolation configurations
on the triangular lattice ${\mathbb T}$, i.e.,
$\text{T}: \Omega \rightarrow \Omega$,
where $\Omega$ is the set of configurations
$\{ -1, +1 \}^{\mathbb T}$.
This family is parametrized by $n \in \{ 0, 1, \ldots \} \cup \{ \infty \}$
representing (discrete) time.
The initial configuration $\sigma^0$ consists of an assignment
of $-1$ or $+1$ to the sites of ${\mathbb T}$.
At times, we will identify the spin variable $\sigma_x$ with the
corresponding site (or hexagon) $x$.
We choose $\sigma^0$ according to a probability measure
$\mu^0$ corresponding to independent identically distributed
$\sigma^0_x$'s with $\mu^0(\sigma^0_x = +1) = \lambda \in[0,1]$.
With the exception of Proposition \ref{perc}, we will set
$\lambda = 1/2$, so that $\mu^0$ is the distribution
corresponding to critical independent site percolation.
We denote by $\mu^n$ the distribution of $\sigma^n$
and write $\mu^n = \tilde{\text{T}}^n \mu^0$,
where $ \tilde{\text{T}}$ is a map in the space of
measures on $\Omega$.
The action of the cellular automaton T can be described as follows:
\begin{equation}
\sigma_x^{n+1} = \left\{ \begin{array}{ll}
                                \sigma_x^n & \mbox{if $x$ has at
least two neighbors $y_1$ and $y_2$ such that
$\sigma_{y_1}^n = \sigma_{y_2}^n = \sigma_x^n$,} \\
        &       \mbox{and $y_1$ and $y_2$ are not neighbors} \\
                                - \sigma_x^n & \mbox{otherwise}
                             \end{array}
                     \right.
\end{equation}
Once the initial percolation configuration $\sigma^0$ is chosen,
the dynamics is completely deterministic, that is, T is a deterministic
cellular automaton with random initial state (see Figs.~\ref{example3}
and~\ref{example4}).

Certain configurations are stable for the
dynamics, in other words, they are absorbing states for the cellular
automaton.
To see this, let us consider a loop in ${\mathbb T}$ expressed as a sequence
of sites $(\zeta_0, \ldots, \zeta_k)$ which are distinct except that
$\zeta_0 = \zeta_k$ and suppose
moreover that $k \geq 6$ and that $\zeta_{i-1}$ and
$\zeta_{i+1}$ are not neighbors; we will call such a sequence an
{\bf m-loop}.
If $\sigma_{\zeta_0} = \sigma_{\zeta_1} = \ldots = \sigma_{\zeta_{k-1}}$
then every site $\zeta_i$ in such an m-loop has two neighbors, $\zeta_{i-1}$ and
$\zeta_{i+1}$, such that $\zeta_{i-1}$ and $\zeta_{i+1}$ are not neighbors
of each other and $\sigma_{\zeta_{i-1}} = \sigma_{\zeta_i} = \sigma_{\zeta_{i+1}}$.
According to the rules of the dynamics, $\sigma_{\zeta_i}$ is
therefore stable, that is, retains the same sign at all future times.
Other stable configurations are ``barbells,'' where a barbell consists
of two disjoint such m-loops connected by a stable m-path.
The stability of certain loops under the action of T will be a key
ingredient in the proof of the main theorem.
A more precise definition of m-paths and m-loops is given in Section
\ref{proofs}.

An important feature of this dynamics is that almost surely every spin
flips only a finite number of times and every local configuration gets
fixated in finite time.
To show this, we introduce a formal Hamiltonian
\begin{equation} \label{energy}
{\cal H}(\sigma) = - \sum_{ \langle x,y \rangle } \sigma_x \sigma_y =
\frac{1}{2} \sum_{x} {\cal H}_x (\sigma),
\end{equation}
where $\sum_{ \langle x,y \rangle }$ denotes the sum over all pairs of 
neighbor sites, each pair counted once, and
\begin{equation} \label{x_energy}
{\cal H}_x (\sigma) =  - \sum_{y \in {\cal N}(x)} \sigma_x \sigma_y,
\end{equation}
where ${\cal N}(x)$ is the set of six (nearest) neighbors of $x$.
We also introduce a ``local energy''
\begin{equation} \label{local_energy}
{\cal H}_{\Lambda}(\sigma) =
- \sum_{ \stackrel{\langle x,y \rangle}{x,y \in \Lambda} } \sigma_x \sigma_y
-  \sum_{z \in \partial \Lambda} \sum_{ \stackrel{x \in \Lambda} {x \in {\cal N}(z)} }
\sigma_x \sigma_z,
\end{equation}
where $\Lambda$ is a subset of ${\mathbb T}$ and $\partial \Lambda$ is
the outer boundary of $\Lambda$, i.e.,
$\{ \zeta \notin \Lambda : x \in {\cal N}(\zeta) \text{ for some } x \in \Lambda \}$.
Notice that although the total energy ${\cal H}(\sigma)$ is almost surely infinite
and is therefore only defined formally, we will only use local energies
of finite subsets of ${\mathbb T}$.

The notion of the energy change caused by a spin flip is somewhat ambiguous
in a cellular automaton  because of the synchronous dynamics and hence
multiple simultaneous spin flips.
Nevertheless, it is easy to show (see the proof of Proposition \ref{fixation})
that each step of the dynamics either lowers or leaves unchanged the energy
-- both locally and globally.
In this sense, our cellular automaton can be considered a zero-temperature
dynamics (see, for example, \cite{nns} and references therein).

\begin{proposition} \label{fixation}
For all values of $\lambda$, almost surely, every spin flips only a
finite number of times.
\end{proposition}

\noindent {\it Proof.} By the translation invariance and ergodicity of the model,
it is enough to prove the claim for the origin.
At time zero, and for all values of $\lambda$, the origin is almost surely
surrounded by an m-loop $\Gamma$ of spins of constant sign.
(For $\lambda = 1/2$, there are infinitely many loops of both signs surrounding
the origin.)
Such an m-loop is stable for the dynamics and its spins retain the same sign
at all times.
Call ${\Lambda}$ the (a.s. finite) region surrounded by $\Gamma$.
The energy ${\cal H}_{\Lambda}$ of the \emph{finite} region $\Lambda$
is bounded below.
Consequently, if we can show that no step of the dynamics ever raises
${\cal H}_{\Lambda}$, it would follow that there can only be a finite number
of steps that strictly lower the energy ${\cal H}_{\Lambda}$.

Call an edge $\eta_{x,y}$ {\bf satisfied} if $\sigma_x = \sigma_y$ and
{\bf unsatisfied} otherwise; then the change in energy
${\cal H}_{\Lambda} (\sigma^{n+1}) - {\cal H}_{\Lambda} (\sigma^n)$
is twice the difference between the number of satisfied edges at time $n$
that become unsatisfied at time $n+1$ and the number of unsatisfied edges
at time $n$ that become satisfied at time $n+1$.
The edges that change from satisfied to unsatisfied or vice-versa are those
between spins that flip and their neighbors that do not, so
\begin{equation} \label{energy_diff}
{\cal H}_{\Lambda} (\sigma^{n+1}) - {\cal H}_{\Lambda} (\sigma^n)
= \sum_{x \in \Lambda \, : \, \sigma_x^{n+1} \neq  \sigma_x^n}
\sum_{ y \in {\cal N}(x) }
(\sigma_x^n
 \sigma_y^n - \sigma_x^{n+1} \sigma_y^{n+1})
= \sum_{x \in \Lambda \,: \, \sigma_x^{n+1} \neq  \sigma_x^n} \Delta_n {\cal H}_x,
\end{equation}
where
\begin{equation}
\Delta_n {\cal H}_x
= {\cal H}_x (\sigma^{n+1}) - {\cal H}_x (\sigma^n) .
\end{equation}
Notice that the only nonzero contributions in the first sum of
(\ref{energy_diff}) come from those sites $y \in \Lambda$ that do not flip
at time $n$.
We want to show that $\Delta_n {\cal H}_x \leq 0$ for all $x \in \Lambda$
and find  some $y$ with $\Delta_n {\cal H}_y < 0$ (assuming there was at
least one spin flip inside $\Lambda$ at time $n$).

Call $D_x^n$ the number of disagreeing neighbors of $x$ at time $n$
and notice that a necessary condition for the spin at site $x$ to flip
at that time is $D^n_x \geq 4$.
Let us first consider the case $D_x^n \geq 5$
and assume, without loss of generality, that $\sigma_x^n = -1$ and
$\sigma_x^{n+1} = +1$.
Then, at time $n$, site $x$ has at least five plus-neighbors,
and at least three of them have plus-spins at time $n+1$ (those
having at time $n$ two plus-neighbors that are not neighbors
of each other).
This implies that the number of edges incident on $x$ that change from
unsatisfied to satisfied is at least three, while the number of edges
that change from satisfied to unsatisfied is one.
Then,
\begin{equation}
\Delta_n {\cal H}_x
= \sum_{y \in {\cal N}(x)}
(\sigma_x^n \sigma_y^n - \sigma_x^{n+1} \sigma_y^{n+1}) < 0.
\end{equation}

\begin{figure}[!ht]
\begin{center}
\includegraphics[height=6cm]{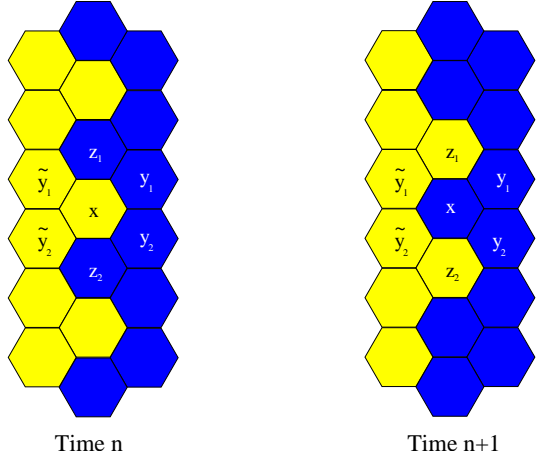}
\caption{Example of a step of the dynamics acting on a local configuration
that leaves the energy of some sites (including $x$) unchanged and decreases
that of other sites (the three spins at the top and the three at the bottom
are stable, as determined also by spins that do not appear in the figure).}
\label{example5}
\end{center}
\end{figure}

Next, we consider the case $D^n_x=4$.
Again, we can assume that $\sigma^n_x = -1 = -\sigma^{n+1}_x$.
In this case, one can have two types of spin flips, one with
${\cal H}_x (\sigma^{n+1}) - {\cal H}_x (\sigma^n) < 0$ and one with
${\cal H}_x (\sigma^{n+1}) - {\cal H}_x (\sigma^n) = 0$.
The second type occurs when,
of the four neighbors of $x$ that are plus at time $n$, two remain
plus at time $n+1$ (we call them $y_1$ and $y_2$)
and the other two flip to minus (we call them $z_1$ and $z_2$), while
the two neighbors that are minus at time $n$ remain minus at time $n+1$
(we call them $\tilde y_1$ and $\tilde y_2$).
In this situation, site $x$ disagrees with four of its six neighbors
both at time $n$ and at time $n+1$ (see Figure \ref{example5}) and therefore
${\cal H}_x (\sigma^{n+1}) - {\cal H}_x (\sigma^n) = 0$.

\noindent The spins at $z_1$ and $z_2$ flip together with the spin at $x$.
Each energy, ${\cal H}_{z_1}$ and ${\cal H}_{z_2}$, is lowered or left
unchanged.
If either is lowered, then the energy ${\cal H}_{\Lambda}$ is lowered.
If neither is lowered, then $z_1$ and $z_2$ must be in the same situation
as $x$ (but with minus and plus interchanged), which requires a configuration
that looks locally like the one in the left part of Figure \ref{example5}
(or one equivalent to it under some lattice symmetry).
However, such a local configuration cannot extend forever; it must be finite
and contained in $\Lambda$.
This implies that we will necessarily find at least one spin $y$ (in fact,
at least two spins) that flip together with the spin at $x$ and such that
$\Delta_n {\cal H}_y < 0$.

Thus, if at time $n$ some site $x \in \Lambda$ flips,
\begin{equation}
{\cal H}_{\Lambda} (\sigma^{n+1}) - {\cal H}_{\Lambda} (\sigma^n)
= \sum_{x \in \Lambda \,: \, \sigma_x^{n+1} \neq  \sigma_x^n} \Delta_n {\cal H}_x < 0.
\end{equation}
It follows that there can only be a finite number of times $n$ at which
spins in $\Lambda$ (in particular, the origin) flip.
This completes the proof. \fbox{} \\

Let us now give two results that are analogous to results proved
in \cite{cns} for a related cellular automaton that will be
discussed below in Section \ref{last}.

\begin{proposition} \label{perc}
If $\lambda > 1/2$ (respectively, $< 1/2$), then for almost every $\sigma^0$,
there is percolation of $+1$ (respectively, $-1$) spins in $\sigma^n$
for any $n \in [0, \infty]$.
\end{proposition}

\noindent {\it Proof.} We only give the proof for $\lambda > 1/2$,
since the case $\lambda < 1/2$ is the same by symmetry.
If $\lambda > 1/2$, since the critical value for independent (Bernoulli)
percolation on the triangular lattice is exactly $1/2$, there is at time
zero percolation of $+1$ spins.
This implies the existence, at time $0$, of doubly-infinite plus-paths
that are stable.
Therefore, there is percolation of $+1$ spins for all $n \geq 0$. \fbox{} \\

We denote by $C_x$ the cluster at $x$ (for a configuration $\sigma$),
i.e., the maximal connected set
$B \in {\mathbb T}$ such that $x \in B$ and $\sigma_y = \sigma_x$
for all $y \in B$.
We write $C_x(n)$ to indicate the cluster at $x$ for $\sigma^n$.
$C_x(n)$ is random and its distribution is denoted $\mu^n$.
$\text{E}_{\mu^n}$ denotes expectation with respect to $\mu^n$.

\begin{proposition} \label{mean_c_s}
For $\lambda = 1/2$, the following two properties are valid. 
\begin{enumerate}
\item For almost every $\sigma^0$, there is no percolation in $\sigma^n$
of either $+1$ or $-1$ spins for any $n \in [0, \infty]$.
\item The mean cluster size in $\sigma^n$ is infinite for any
$n \in [0, \infty]$ \emph{:} for any $x \in {\mathbb T}$,
\begin{equation} \label{inf_mean}
\emph{E}_{\mu^n} (|C_x|) = \emph{E}_{\mu^0} (|C_x(n)|) = \infty.
\end{equation}
\end{enumerate}
\end{proposition}

\noindent {\it Proof.}
To prove the first claim, notice that at time zero the origin is
almost surely surrounded by both a plus and a minus m-loop.
Those loops are stable and prevent the cluster at the origin,
be it a plus or a minus cluster, from percolating at all subsequent
times.
Therefore the probability that the origin belongs to an infinite cluster
is zero for all $n \geq 0$.

To prove the second claim, we note that because of the absence of
percolation for either sign, it follows from a theorem of Russo~\cite{russo}
applied to the triangular lattice, that the mean cluster size of both
plus and minus clusters diverges. \fbox{} \\

Before stating our main theorem, we need some more definitions
to formulate the continuum scaling limit.
We adopt the approach of~\cite{ab} (see also~\cite{abnw}).

\subsection{Compactification of ${\mathbb R}^2$} \label{comp}

The scaling limit $\delta \to 0$ can be taken by focusing on fixed
finite regions, $\Lambda \subset {\mathbb R}^2$, or by treating the
whole ${\mathbb R}^2$.
The second option is more convenient, because it avoids technical
issues that arise near the boundary of $\Lambda$.

A convenient way of dealing with the whole ${\mathbb R}^2$
is to replace the Euclidean metric with a distance function
$\text{d}(\cdot,\cdot)$ defined on ${\mathbb R}^2 \times {\mathbb R}^2$
by
\begin{equation} \label{metric}
\text{d}(u,v) = \inf_{\phi} \int (1 + |\phi|^2)^{-1} \, ds ,
\end{equation}
where the infimum is over all smooth curves $\phi(s)$
joining $u$ with $v$, parametrized by arc-length $s$,
and $|\cdot|$ denotes the Euclidean norm.
This metric is equivalent to the Euclidean metric in bounded regions,
but it has the advantage of making ${\mathbb R}^2$ precompact.
Adding a single point at infinity yields the compact space
$\dot{\mathbb R}^2$ which is isometric, via stereographic projection,
to the two-dimensional sphere.

\subsection{The space of curves} \label{curves}

Denote by ${\cal S}$ the complete, separable metric space of continuous
curves in $\dot{\mathbb R}^2$ with a distance $\text{D}(\cdot,\cdot)$
based on the metric defined by eq.~(\ref{metric}) as follows.
Curves are regarded as equivalence classes of continuous functions $\gamma(t)$
from the unit interval to $\dot{\mathbb R}^2$, modulo monotonic
reparametrizations.
${\cal C}$ will represent a particular curve and $\gamma(t), t \in [0,1]$, a
particular parametrization of ${\cal C}$, while ${\cal F}$ will represent a
set of curves.
The distance D between two curves, ${\cal C}_1$ and ${\cal C}_2$, is defined by
\begin{equation}
\text{D} ({\cal C}_1,{\cal C}_2) \equiv \inf_{f_1,f_2} \sup_{t \in [0,1]}
\text{d} (\gamma_1 (f_1(t)), \gamma_2 (f_2(t))),
\end{equation}
where $\gamma_1$ and $\gamma_2$ are particular parametrizations of
${\cal C}_1$ and ${\cal C}_2$, and the infimum is over the set of all monotone
(increasing or decreasing) continuous functions from the unit interval
onto itself.
The distance between two closed sets of curves is defined by the
induced Hausdorff metric as follows:
\begin{equation} \label{hausdorff}
\text{dist}({\cal F},{\cal F}') \leq \varepsilon
\Leftrightarrow \forall \, {\cal C} \in {\cal F}, \, \exists \,
{\cal C}' \in {\cal F}' \text{ with }
\text{D} ({\cal C},{\cal C}') \leq \varepsilon,
\text{ and vice-versa.}
\end{equation}

For each fixed $\delta>0$, the random curves that we consider are polygonal
paths in the hexagonal lattice $\delta {\mathbb H}$, dual of $\delta {\mathbb T}$,
consisting of connected portions of the boundaries between plus and minus
clusters in $\delta {\mathbb T}$.
A subscript $\delta$ may then be added to indicate that the
curves correspond to a model with a ``short distance cutoff.''
The probability measure $\mu^n_{\delta}$
denotes the distribution of the random set of curves ${\cal F}^n_{\delta}$
consisting of the polygonal paths on
$\delta {\mathbb H}$ generated by $\text{T}^n$ acting on $\sigma^0$.


\subsection{Main results} \label{main}

Since the cellular automaton is deterministic, all the percolation models
$\sigma^n$ for $n \geq 0$ are automatically coupled, once they are all
constructed on the single probability space $( \Omega, \Sigma, \mu^0)$
on which $\sigma^0$ is defined.
The following theorem is valid for any $\lambda$, but $\lambda = 1/2$
is the only interesting case, therefore we restrict attention to it.

\begin{theorem} \label{thm}
For $\lambda = 1/2$, the Hausdorff distance between the system of random
curves ${\cal F}^0_{\delta}$ at time $0$ and the corresponding system of
curves ${\cal F}^n_{\delta}$ at time $n$ goes to zero almost surely
as $\delta \to 0$, for each $n \in [1,\infty]$;
i.e., for $\mu^0$-almost every $\sigma^0$,
\begin{equation} \label{convergence}
\lim_{\delta \to 0} \text{\emph{dist}}
({\cal F}^0_{\delta},{\cal F}^n_{\delta}) = 0,
\text{ for any } n \in [1, \infty].
\end{equation}
\end{theorem}

\noindent The main application of the theorem is that the scaling limits
of our family ${\cal F}_1^n$ of percolation models, if they exist, must be
the same for all $n \in [0, \infty]$:

\begin{corollary} \label{cor}
Suppose that critical site percolation on the triangular lattice
has a unique scaling limit in the Aizenman-Burchard sense \cite{ab},
i.e., ${\cal F}_{\delta}^0$ converges in distribution as $\delta \to 0$
(for $\lambda = 1/2$).
Then, for every $n \in [1, \infty]$, ${\cal F}_{\delta}^n$ converges
in distribution to the same limit.
\end{corollary}


\section{Proofs} \label{proofs}

In this section, we give the proofs of the main results.
We start by reminding the reader of the definitions of m-path and m-loop
and by giving some new definitions and two lemmas which will be used later.

\begin{definition}
A {\bf path} $\Gamma$ between $x$ and $y$ in $\delta {\mathbb T}$
(embedded in ${\mathbb R}^2$) is an ordered sequence of sites
$(\zeta_0=x, \ldots, \zeta_k=y)$ with $\zeta_i \neq \zeta_j$ for
$i \neq j$ and $\zeta_{i+1} \in {\cal N}(\zeta_i)$.
A {\bf loop} is a sequence $(\zeta_0, \ldots, \zeta_{k+1}=\zeta_0)$
with $k \geq 1$ such that $(\zeta_0, \ldots, \zeta_k)$ and
$(\zeta_1, \ldots, \zeta_{k+1})$ are paths.

We call a path  $(\zeta_0, \ldots, \zeta_k)$ (resp., a loop
$(\zeta_0, \ldots, \zeta_{k+1} = \zeta_0)$) an {\bf m-path}
(resp., an {\bf m-loop}) if $\zeta_{i-1}$ and $\zeta_{i+1}$ are not
neighbors for $i = 1, \ldots, k-1$ (resp., for $i = 1, \ldots, k+1$,
where $\zeta_{k+2} = \zeta_1$).
\end{definition}

\noindent Notice that for every path $\Gamma$ between $x$ and $y$,
there always exists at least one m-path $\tilde\Gamma$ between the
same sites.

\begin{definition}
A {\bf boundary path} ({\bf b-path}) $\Gamma^*$ is an ordered sequence
$(\eta^*_0 = \eta^*_{\zeta_0, \xi_0}, \ldots,
\eta^*_k = \eta^*_{\zeta_k, \xi_k})$ of distinct dual edges such
that either $\zeta_{i+1} = \zeta_i$ and $\xi_{i+1} \in {\cal N}(\xi_i)$
or $\xi_{i+1} = \xi_i$ and $\zeta_{i+1} \in {\cal N}(\zeta_i)$, and
$\sigma_{\zeta_i} = - \sigma_{\xi_i}$ for all $i = 0, \ldots, k$.

We call a maximal boundary path simply a {\bf boundary}; it can either
be a doubly-infinite path or a finite loop.
\end{definition}

b-paths represent the (random) curves introduced in
the previous section for a fixed value of the short distance cutoff
$\delta$.
A collection of such curves in the compact space $\dot{\mathbb R}^2$
is indicated by ${\cal F}_{\delta}$.

b-paths $\Gamma^*$ are parametrized by functions $\gamma (t)$,
with $t \in [0,1]$.
When we write that, for $t$ in some interval $[t_1,t_2]$ (the interval
could as well be open or half-open),
$\gamma (t) \in \eta^*_{x,y}$ 
we mean that the parametrization $\gamma (t)$ for $t$ between
$t_1$ and $t_2$ is irrelevant and can be chosen, for example,
so that $|\frac{d\gamma (t)}{dt}|$ is constant for $t \in [t_1,t_2]$.
The notation $\Gamma^* (u, v)$, where $u$ and $v$ can be dual sites
or generic points of ${\mathbb R}^2 \cap \Gamma^*$,
stands for the portion of $\Gamma^*$ between $u$ and $v$.

\begin{definition}
A {\bf stable loop} $l$ (for some $\sigma$) is an m-loop
$(\zeta_0=x, \ldots, \zeta_k=x)$
such that $\sigma_{\zeta_0} = \sigma_{\zeta_1} = \ldots = \sigma_{\zeta_k}$.
\end{definition}

\begin{definition}
We say that a dual edge $\eta^*_{x, y}$ is {\bf stable} if $x$ and $y$
belong to stable loops of opposite sign.
\end{definition}

\begin{lemma} \label{lemma1}
In general, a constant sign m-path $(\zeta_0, \ldots, \zeta_k)$
is ``fixated'' (i.e., retains that same sign in $\sigma^n$ for all
$0 \leq n \leq \infty$) if $\zeta_0$ and $\zeta_k$ are fixated
\end{lemma}

\noindent {\it Proof}. The claim of the lemma is straightforward.

\begin{lemma} \label{lemma2}
For $\lambda = 1/2$,
there is a one to one mapping from boundaries in ${\cal F}^{n+1}_{\delta}$
to ``parent'' boundaries in ${\cal F}^n_{\delta}$.
\end{lemma}

\noindent {\it Proof}. Let $\Gamma_n^* \in {\cal F}^n_{\delta}$
be a boundary at time $n$ and $\text{int}(\Gamma_n^*)$ be the set
of sites of $\mathbb T$ (i.e., hexagons) ``surrounded'' by $\Gamma^*_n$.
Let $C$ the (unique) constant sign cluster contained in $\text{int}(\Gamma_n^*)$
that has sites next to $\Gamma_n^*$.
Suppose that sites $x$ and $y$ in $C$ do not change sign at time $n+1$,
then, at that time, they must belong to the same cluster $C'$.
(This means that a cluster cannot split in two or more pieces under the effect
of the dynamics.)
The reason is that $x$ and $y$ are connected by an m-path $\Gamma$ at time
$n$ because they are in the same cluster, and since they have not flipped
between time $n$ and time $n+1$, Lemma~\ref{lemma1} (or more accurately,
the single time step analogue of Lemma~\ref{lemma1}) implies that the
sites in $\Gamma$ have not flipped either, so that $x$ and $y$ are still
connected at time $n+1$ and therefore belong to the same cluster.

If at least one site $x$ in $C$ does not change sign at time $n+1$,
we say that $C$ has survived (and evolved into a new cluster $C'$ that contains
$x$ at time $n+1$).
From what we just said, the sites of $C$ that retain the same sign form
a unique cluster $C'$.
We call $C$ the parent cluster of $C'$ and the external boundary
$\Gamma_n^* \in {\cal F}^n_{\delta}$ of $C$ the parent of the external boundary
$\Gamma_{n+1}^* \in {\cal F}^{n+1}_{\delta}$ of $C'$ (note that $C'$ is a.s. finite
-- see Proposition~\ref{mean_c_s}).
To prove that there is a one to one mapping between boundaries in ${\cal F}^{n+1}_{\delta}$
and those in ${\cal F}^n_{\delta}$, it remains to show that the parent boundary
of each element of ${\cal F}^{n+1}_{\delta}$ is unique.
(This means that clusters cannot merge under the effect of the dynamics.)

If this were not the case, then a cluster $C'$ could have two or more
distinct parent clusters, $C_1, C_2, \ldots$ .
Notice that each of the parent clusters at time $n$ is surrounded by a
constant sign m-loop $\Gamma_i$, $i = 1, 2, \ldots$, which is stable.
Suppose that $x_1 \in C_1$ and $x_2 \in C_2$ retain their sign at time
$n+1$.
Without loss of generality, we shall assume that this sign is plus.
Assuming that $C_1$ and $C_2$ are both parents of $C'$, $x_1$ and $x_2$
should both belong to $C'$ at time $n+1$ and therefore be connected
by a plus-path.
But this contradicts the fact that, at time $n$, $x_1 \in \text{int}(\Gamma_1)$
and $x_2 \in \text{int}(\Gamma_2)$, that is, they belong to the interiors
of two disjoint stable minus-loops. \fbox{} \\

Notice that boundaries cannot be ``created,'' but can ``disappear,''
as complete clusters are ``eaten'' by the dynamics.

\subsection{Proof of Theorem \ref{thm}} \label{proof-thm}

Let us start, for simplicity, with the case of a single b-path.
For any $\varepsilon > 0$, given a b-path
$\Gamma^*_0$ at time $0$, parametrized by some $\gamma(t)$,
we will find a path $\Gamma^*_n$ at time $n$
with parametrization $\gamma'(t)$ such that, for $\delta$ small enough,
\begin{equation} \label{small-distance}
\sup_{t \in [0,1]} \text{d}(\gamma(t), \gamma'(t)) \leq \varepsilon + 2 \delta
\end{equation}
and vice-versa (i.e., given $\Gamma^*_n$ and $\gamma'$, we need can
find $\Gamma^*_0$ and $\gamma$ so that (\ref{small-distance}), in which
the dependence on the scale factor $\delta$ has been suppressed, is valid).
Later we will require that this holds simultaneously for all the
curves in ${\cal F}^0_{\delta}$ and ${\cal F}^n_{\delta}$,
as required by eq. (\ref{convergence}).

Let $B^1(R)$ be the ball of radius $R$, $B^1(R)=\{u \in {\mathbb R}^2 :
|u| \leq R\}$ in the Euclidean metric, and
$B^2(R)=\{u \in {\mathbb R}^2 : \text{d}(u) \leq R\}$ the ball of radius $R$
in the metric~(\ref{metric}).
For a given $\varepsilon > 0$, we divide ${\mathbb R}^2$ into two regions:
$B^1(6/\varepsilon)$ and ${\mathbb R}^2 \setminus B^1(6/\varepsilon)$.
We start by showing that, thanks to the choice of the metric~(\ref{metric}),
one only has to worry about curves (or polygonal paths) that
intersect $B^1(6/\varepsilon)$.
In fact, the distance between any two points $u,v \in
\dot{\mathbb R}^2 \setminus B^1(6/\varepsilon)$ satisfies the
following bound
\begin{equation}
\text{d}(u,v) \leq \text{d}(u,\infty) + \text{d}(v,\infty) \leq
2 \int_0^{\infty} [1 + (s + 6 / \varepsilon)^2]^{-1} \, ds < \varepsilon/3.
\end{equation}
Thus, given any curve in ${\cal F}^0_{\delta}$ contained completely in
$\dot{\mathbb R}^2 \setminus B^1(6/\varepsilon)$, it can be approximated
by any curve in ${\cal F}^n_{\delta}$ also contained in
$\dot{\mathbb R}^2 \setminus B^1(6/\varepsilon)$, and viceversa.
The existence of such curves in ${\cal F}^0_{\delta}$ is clearly not a problem,
since the region $\dot{\mathbb R}^2 \setminus B^1(6/\varepsilon)$
contains an infinite subset of $\delta \mathbb T$ and therefore there
is zero probability that it doesn't contain any b-path at time zero.
There is also zero probability that it contains no stable b-path at time zero,
but any such b-path also belongs to ${\cal F}^n_{\delta}$.

Before we can proceed, we need the following lemma, which is a
consequence of the the fact that at time zero we are dealing with
a Bernoulli product measure.
In this lemma (and elsewhere), the {\bf diameter} $\text{diam}(\cdot)$
of a subset of ${\mathbb R}^2$ is defined by using the Euclidean metric.

\begin{lemma} \label{exp-bound}
Let $\eta_0^*$ be any (deterministic) dual edge; then for
some constant $c>0$,

\begin{equation}
\mu^0 \left( \exists \, \Gamma^* \ni \eta_0^* :
\text{\emph{diam}} (\Gamma^*) \geq M \text{ and } \Gamma^*
\text{ does not contain at least one stable edge} \right) 
\leq e^{- c M}.
\end{equation}
\end{lemma}

\begin{figure}[!ht]
\begin{center}
\includegraphics[width=8cm]{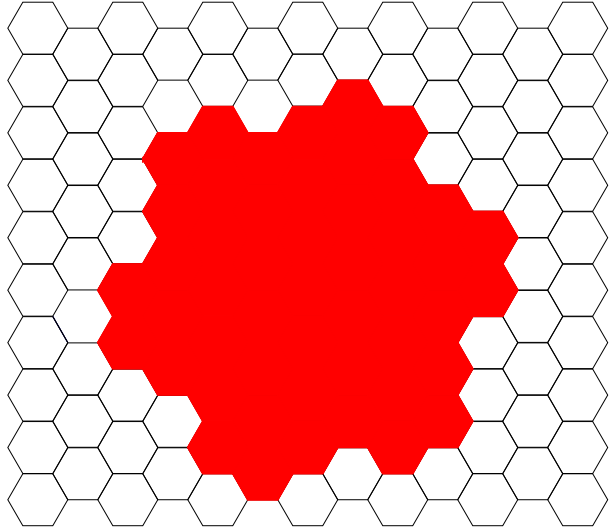}
\caption{Elementary cell for the partition of ${\mathbb H}$ used
in Lemma \ref{exp-bound}. Notice that the cell is made out of seven
smaller cells, each of them formed by seven hexagons.
}\label{fig-lemma1}
\end{center}
\end{figure}

\noindent {\it Proof}.
To prove the lemma, we partition the hexagonal lattice into
regions $Q_i$ as in Figure~\ref{fig-lemma1}.
We then do an algorithmic construction of $\Gamma^*$, starting
from $\eta_0^*$, as a percolation
exploration process, but with the additional rule that, when the
exploration process hits the boundary of any $Q_i$ for the first time,
all the hexagons in $Q_i$ are checked next
(according to some deterministic order).
From every entrance point of $Q_i$, there is a choice of the
values of the spins of the outermost layer of $Q_i$ that forces
$\Gamma^*$ to enter $Q_i$.
Therefore, when $\Gamma^*$ hits $Q_i$, it always has a positive
probability of entering the region.

We call $F_i$ the event that a (dual) stable edge is found inside $Q_i$,
belonging to $\Gamma^*$.
It is easy to see that such an event has positive probability, bounded
away from zero by a constant that does not depend on how the exploration
process enters the region $Q_i$.
In fact, from any entrance point, there is clearly a choice of the values
of the spins in $Q_i$ that forces $\Gamma^*$ to cut $Q_i$ in two symmetric
parts, containing spins of opposite sign.

Now, if $\text{diam}(\Gamma^*) \geq M$, then $\Gamma^*$ must clearly
visit at least $O(M)$ different regions $Q_i$.
The conclusion of the proof should now be clear (see, for example,
\cite{fn}). \fbox{} \\

With this lemma, we can now proceed to the proof of the theorem.
As explained before, we restrict attention to paths that intersect
$B^1(6/\varepsilon)$.
Given a b-path $\Gamma^*_0 
= (\eta^*_0 = \eta^*_{\zeta_0,\xi_0}, \ldots, \eta^*_k = \eta^*_{\zeta_k,\xi_k})$
in ${\cal F}^0_{\delta}$ with parametrization $\gamma(t)$,
call $u_0$ the point in ${\mathbb R}^2$ where the first (dual)
edge $\eta^*_0$ begins.
The following algorithmic construction produces a sequence
$u_0, \ldots, u_N$ of points in $\Gamma^*_0$.
\begin{enumerate}
\item Start with $u_0$.
\item Once $u_0, \ldots, u_i$ have been constructed,
if $u_i \in B^1(6/\varepsilon)$
take the ball $B^1_{u_i}(\varepsilon / 3)$ centered at $u_i$ and
of radius $\varepsilon / 3$ and let $u_{i+1}$ be the first intersection
of $\Gamma^*_0 \setminus \Gamma^*_0 (u_0, u_i)$ with
$\partial B^1_{u_i}(\varepsilon / 3)$,
if $u_i \notin B^1(6/\varepsilon)$
take the ball $B^2_{u_i}(\varepsilon / 3)$ centered at $u_i$ and
of radius $\varepsilon / 3$ and let $u_{i+1}$ be the first intersection
of $\Gamma^*_0 \setminus \Gamma^*_0 (u_0, u_i)$ with
$\partial B^2_{u_i}(\varepsilon / 3)$.
\item Terminate when there is no next $u_i$.
\end{enumerate}

During the construction of the sequence $u_0, \ldots, u_N$,
$\Gamma^*_0$ is split in $N+1$ pieces, the first $N$ having
diameter at least $\varepsilon / 3$.
The construction also produces a sequence of balls
$B^{j_0}_{u_0}(\varepsilon / 3), \ldots, B^{j_N}_{u_N}(\varepsilon / 3)$,
with $j_i = 1$ or $2$.
Notice that no two successive $u_i$'s can lie outside of
$B^1(6 / \varepsilon)$.
In fact, if for some $i$, $u_i$ lies outside of $B^1(6 / \varepsilon)$,
$u_{i+1}$ belongs to $\partial B^2_{u_i}(\varepsilon / 3)$,
which is contained inside $B^1(6 / \varepsilon)$, due to the choice
of the metric.
Each $u_i$ contained in $B^1(6 / \varepsilon)$ lies on an edge of
$\delta {\mathbb H}$, but no more than one $u_i$ can lie on the
same edge since $\Gamma^*_0$ is self-avoiding and cannot use
the same edge or site more than once.
Also, the number of $u_i$'s lying outside of $B^1(6 / \varepsilon)$
cannot be larger than (one more than) the number of the $u_i$'s lying
inside $B^1(6 / \varepsilon)$.
Therefore, $N \leq const \times (\varepsilon \delta)^{-2}$.

For any two successive balls, $B^{j_i}_{u_i}(\varepsilon / 3)$ and
$B^{j_{i+1}}_{u_{i+1}}(\varepsilon / 3)$, let $O_i =
B^{j_i}_{u_i}(\varepsilon / 3 + \delta)
\cup B^{j_{i+1}}_{u_{i+1}}(\varepsilon / 3 + \delta)$.
Now assume that there exists a sequence
$\bar \eta^*_0 ,\ldots \bar \eta^*_{N-1}$ of stable (dual) edges
of $\Gamma^*_0$, with $\bar \eta^*_i$ contained in
$\Gamma^*_0 (u_i, u_{i+1})$.
$\Gamma^*_0 (\bar \eta^*_i, \bar \eta^*_{i+1})$ is contained in
$O_i$ (for fixed $\varepsilon$ and small enough $\delta$).
Also contained in $O_i$ are the two paths (i.e., paths in
$\delta {\mathbb T}$ which may be thought as sequences of hexagons)
whose hexagons are next to $\Gamma^*_0 (\bar \eta^*_i, \bar \eta^*_{i+1})$,
one on each side.
From those two paths one can extract two subsets that are m-paths 
which are stable since the first and last hexagons of each one of them
is stable (such hexagons must be stable since they are next to
$\bar\eta^*_i$ and $\bar\eta^*_{i+1}$).
The two m-paths so constructed constitute a ``barrier'' that limits
the movements of the boundary, so that
$\Gamma^*_n (\bar \eta^*_i, \bar \eta^*_{i+1})$
is in fact confined to lie within those two m-paths and thus within $O_i$.
To parametrize $\Gamma^*_n$, we use any parametrization
$\gamma'(t)$ such that $\gamma'(t) = \gamma(t)$ whenever
$\gamma(t) \in \bar \eta^*_i$.
Using this parametrization and the previous fact, it is clear that
the distance between $\Gamma^*_0 (\bar \eta^*_i, \bar \eta^*_{i+1})$
and $\Gamma^*_n (\bar \eta^*_i, \bar \eta^*_{i+1})$ does not
exceed $\varepsilon + 2 \delta$, the diameter of $O_i$.
Therefore, conditioning on the existence of the above sequence
$\bar \eta^*_0, \ldots, \bar \eta^*_{N-1}$ of stable (dual) edges
of $\Gamma^*_0$, we can conclude that
\begin{equation} \label{small-distance1}
\sup_{t \in [0,1]} \text{d}(\gamma(t), \gamma'(t)) \leq \varepsilon + 2 \delta.
\end{equation}

It remains to prove the existence of the sequence
$\bar \eta^*_0, \ldots, \bar \eta^*_{N-1}$ of stable (dual) edges.
To do that, let us call $A_i$ the event that $\Gamma^*_0 (u_i, u_{i+1})$
does \emph{not} contain at least one stable edge, and let
$A = \cup_{i=0}^{N-1} A_i$ be the event that at least one of
the first $N$ pieces of $\Gamma^*_0$ does not have any
stable edge.
Then, considering that the total number of edges contained in
$B^1(6 / \varepsilon)$ is bounded by $const \times (\varepsilon \delta)^{-2}$
and using Lemma \ref{exp-bound}, we have
\begin{equation} \label{small-prob}
\mu^0_{\delta} (A) 
\leq (\varepsilon \delta)^{-2} e^{- c' (\varepsilon / \delta)}
\end{equation}
for some $c' > 0$.
Equation (\ref{small-prob}) means that the probability of
not finding at least one stable edge in each of the first
$N$ pieces of $\Gamma^*_0$ is very small and goes to $0$,
for fixed $\varepsilon$, as $\delta \to 0$.
This is enough to conclude that, with high probability
(going to $1$ as $\delta \to 0$), equation (\ref{small-distance1})
holds.

This proves one direction of the claim, in the case of a single curve.
To obtain the other direction, notice that a large
b-path at time $n$ is part of a complete boundary $\Gamma^*_n$ which
must come from a line of ``ancestors'' (see Lemma~\ref{lemma2}) that
starts with some $\Gamma^*_0$ at time $0$.
Therefore, one can apply the above arguments to $\Gamma^*_0$,
provided that the latter is large enough.
Although the proof of this last fact is very simple, it is convenient
to state it as a separate lemma.

\begin{lemma} \label{ancestor}
Set $\delta = 1$ for simplicity; then
for any boundary $\Gamma^*_n$ at time $n$, there is an ancestor
$\Gamma^*_0$, with 
$\text{diam}(\Gamma^*_0) \geq \text{diam}(\Gamma^*_n) - 1$.
\end{lemma}

\noindent {\it Proof.}
The existence of an ancestor $\Gamma^*_0$ comes from Lemma \ref{lemma2},
so we just have to show that
$\text{diam}(\Gamma^*_0) \geq \text{diam}(\Gamma^*_n) - 1$.
$\Gamma^*_0$ is surrounded by a connected set of hexagons that touch
$\Gamma^*_0$ and whose spins are all the same (this is the ``external
boundary'' of the set of hexagons that are in the interior of $\Gamma^*_0$).
From this set, one can extract an m-path of constant sign whose diameter
is bounded above by $\text{diam}(\Gamma^*_0) + 1$.
Since such a constant sign m-path is stable for the dynamics,
$\Gamma^*_n$ must lie within its interior.
This concludes the proof. \fbox{} \\

At this point, we need to show that the above argument can be
repeated and the construction done simultaneously for all curves
in ${\cal F}_{\delta}^0$ and ${\cal F}_{\delta}^n$ (for each $n$).
First of all notice that, for a fixed $\varepsilon$, any b-path
$\Gamma^*$ of diameter less than $\varepsilon / 2$ can be approximated
by a closest stable edge, provided that one is found within
the ball of radius $\varepsilon / 2$ that contains the $\Gamma^*$,
with the probability of this last event clearly going to $1$ as
$\delta \to 0$, when we restrict attention to $B^1(6 / \varepsilon)$.
For a b-path outside $B^1(6 / \varepsilon)$, we already noticed
that it can be approximated by any other b-path also outside
$B^1(6 / \varepsilon)$.
As for the remaining b-paths, notice that the total number of
boundaries that intersect the ball $B^1(6 / \varepsilon)$
cannot exceed $const \times (\varepsilon \delta)^{-2}$
(in fact, the total number of pieces in which the boundaries that
intersect $B^1(6 / \varepsilon)$ are divided cannot exceed
$const \times (\varepsilon \delta)^{-2}$).
So, we can carry out the above construction simultaneously for
all the boundaries that touch $B^1(6 / \varepsilon)$, having
to deal with at most $const \times (\varepsilon \delta)^{-2}$
segments of b-paths of diameter of order at least $\varepsilon$.
Therefore, letting
$Y^n_{\delta} = \text{dist}({\cal F}^0_{\delta},{\cal F}^n_{\delta})$,
we can apply once again Lemma \ref{exp-bound} and conclude that
\begin{equation} \label{small-prob1}
\mu^0(Y^n_{\delta} > \varepsilon)
\leq (\varepsilon \delta)^{-2} e^{- c'' (\varepsilon / \delta)}.
\end{equation}

To show that $Y^n_{\delta} \to 0$ as $\delta \to 0$ $\mu^0$-almost surely
and thus conclude the proof,
it suffices to show that, $\forall \varepsilon>0$,
$\mu^0(\limsup_{\delta \to 0} Y^n_{\delta} > \varepsilon) = 0$.
To that end, first take a sequence $\delta_k = 1/2^k$ and
notice that
\begin{equation} \label{bc}
\sum_{k=0}^{\infty} \mu^0(Y^n_{\delta_k} > \varepsilon)
\leq \sum_{k=0}^{\infty} \frac{4^k}{\varepsilon^2} e^{- c'' 2^k \varepsilon}
< \infty,
\end{equation}
where we have made use of (\ref{small-prob1}).
Equation (\ref{bc}) implies that we can apply the Borel-Cantelli lemma and
deduce that
$\mu^0(\limsup_{k \to \infty} Y^n_{\delta_k} > \varepsilon) = 0$,
$\forall \varepsilon>0$.
In order to handle the values of $\delta$ not in the sequence $\delta_k$,
that is for those $\delta$ such that $\delta_{k+1} < \delta < \delta_k$ for some $k$,
we use the following double bound, valid for any $0< \alpha < 1$,
\begin{equation} \label{db}
\alpha \text{d}(u,v) \leq \text{d}(\alpha u, \alpha v) 
\leq \frac{1}{\alpha} \text{d}(u,v),
\end{equation}
which implies that $\alpha Y^n_{\delta_k} \leq Y^n_{\alpha \delta_k}
\leq \frac{1}{\alpha} Y^n_{\delta_k}$.
The two bounds in equation (\ref{db}) come from writing
$\text{d}(\alpha u, \alpha v)$ as
$\text{d}(\alpha u, \alpha v) = \inf_{\phi'} \int (1 + |\phi'|^2)^{-1} \, ds'
= \alpha \inf_{\phi} \int (1 + \alpha^2 |\phi|^2)^{-1} \, ds $,
where $\phi'(s')$ are smooth curves joining $\alpha u$ with $\alpha v$,
while $\phi(s)$ are smooth curves joining $u$ with $v$.

The proof of the theorem is now complete. \fbox{} 

\subsection{Proof of Corollary \ref{cor}}

The Corollary is an immediate consequence of Theorem \ref{thm} and of the
following general fact, of which we include the proof for completeness.

\begin{lemma} \label{general-fact}
If $\{ X_{\delta} \}, \{ Y_{\delta} \}$ (for $\delta > 0$), and $X$ are
random variables taking values in a complete,
separable metric space $S$ (whose $\sigma$-algebra is the Borel algebra)
with $\{ X_{\delta} \}$ and $\{ Y_{\delta} \}$ all defined on the same
probability space, then if $X_{\delta}$ converges in distribution to $X$ and
the metric distance between $X_{\delta}$ and $Y_{\delta}$ tends to zero
almost surely  as $\delta \to 0$, $Y_{\delta}$ also converges in distribution
to $X$.
\end{lemma}

\noindent {\it Proof.} Since $X_{\delta}$ converges to $X$ in distribution,
the family $\{ X_{\delta} \}$ is relatively compact and therefore tight by
an application of Prohorov's Theorem (using the fact that $S$ is a
complete, separable metric space -- see, e.g., \cite{billingsley}).
Then, for any bounded, continuous, real function $f$ on $S$,
and for any $\varepsilon > 0$, there exists a compact set $K$
such that $\int |f(X_{\delta})| I_{\{ X_{\delta} \notin K \}} dP < \varepsilon$
and $\int |f(Y_{\delta})| I_{\{ X_{\delta} \notin K \}} dP < \varepsilon$ for all
$\delta$, where $I_{\{ \cdot \}}$ is the indicator function and $P$ the
probability measure of the probability space of $\{ X_{\delta} \}$ and
$\{ Y_{\delta} \}$.
Thus, for small enough $\delta$,
\begin{equation}
|\int f(X_{\delta}) dP - \int f(Y_{\delta}) dP| <
\int |f(X_{\delta}) - f(Y_{\delta})| I_{\{ X_{\delta} \in K \}} dP + 2 \varepsilon
< 3 \varepsilon,
\end{equation}
where in the last inequality we use the uniform continuity of $f$ when
restricted to the compact set $K$ and the fact that the metric distance
between $X_{\delta}$ and $Y_{\delta}$ goes to $0$ as $\delta \to 0$. \fbox{} \\

To conclude the proof of the Corollary, it is enough to apply
Lemma \ref{general-fact} to $\{ \mu^0_{\delta} \}_{\delta}$,
$\{ \mu^n_{\delta} \}_{\delta}$, $\mu_{sl}$ (or, to be more precise,
to the random variables of which those are the distributions),
for each $n \in [1, \infty]$, where $\mu_{sl}$ is the unique scaling
limit of critical site percolation on the triangular lattice. \fbox{}

\section{Dependent site percolation models on the hexagonal
and triangular lattice} \label{last}

The model that we have presented and discussed in Section~\ref{def+res}
has been chosen as a sort of paradigm, but is not the only one for which
such results can be proved.
In fact, it is not the original model for which such results
were obtained.

In this section we describe some percolation models on the \emph{hexagonal}
lattice and prove that they have the same scaling limit as critical
(independent) site percolation on the \emph{triangular} lattice.
None are independent percolation models, but nonetheless,
they represent explicit examples of critical percolation
models on different lattices with the same scaling limit.
Besides, the construction of the models on the hexagonal lattice
can be seen as a simple and natural way of producing percolation
models for which all the sites of the external (site) boundary of
any constant sign cluster $C$ belong to a unique cluster $C'$ of
opposite sign.
In other words, this implies that the the boundaries between clusters
of opposite sign form a nested collection of loops,
a property that site percolation on the triangular lattice possesses
automatically because of the \emph{self-matching} property of ${\mathbb T}$
(which is crucial in Smirnov's proof of the existence and uniqueness
of the scaling limit of crossing probabilities and Cardy's formula).

\subsection{The models}

The percolation models that we briefly describe here can be
constructed by means of a natural \emph{zero-temperature Glauber dynamics}
which is the zero-temperature case of \emph{Domany's stochastic
Ising ferromagnet} on the hexagonal lattice~\cite{domany}.
The cellular automaton (i.e., Domany's stochastic Ising
ferromagnet at zero temperature) that gives rise to those
percolation models can also be realized on the triangular
lattice with flips when a site disagrees with six,
five and sometimes four of its six neighbors.
The initial state $\sigma^0$ consists of an assignment of $-1$ or $+1$
with equal probability to each site of the hexagonal or
triangular lattice (depending on which version of the cellular
automaton we are referring to).
In the first version, ${\mathbb H}$,
as a bipartite graph, is partitioned into two subsets 
${\cal A}$ and ${\cal B}$ which are alternately updated so that each
$\sigma_x$ is forced to agree with a majority of its three neighbors
(which are in the other subset). 
In the second version, all sites are updated simultaneously according
to a rule based on a deterministic pairing of the six neighbors of
every site into three pairs.
The rule is that $\sigma_x$ flips if and only if it disagrees 
with both sites in two or more
of its three neighbor pairs; thus there is (resp., is not)
a flip if the number $D_x$ of disagreeing neighbors 
is $\ge 5$ (resp., $\le 3$) and there is also a flip for some
cases of $D_x = 4$.
These percolation models on ${\mathbb H}$ and ${\mathbb T}$
are investigated in~\cite{cns1},
where Cardy's formula for rectangular crossing probabilities
is proved to hold in the scaling limit.

\begin{figure}[!ht]
\begin{center}
\includegraphics[width=8cm]{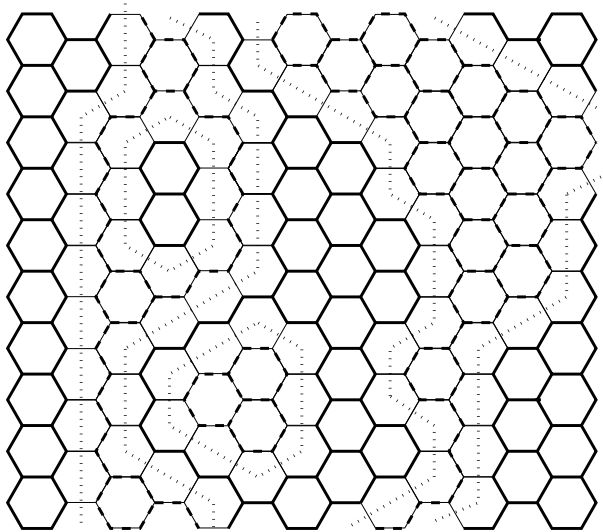}
\caption{Example of a (local) stable configuration for the zero-temperature Domany dynamics.
Heavy lines on edges of $\mathbb H$ connect, say, plus spins, while heavy broken lines connect
minus spins. The dotted lines drawn on edges of the dual lattice are the perpendicular bisectors
of unsatisfied edges, indicating the boundaries between plus and minus clusters.
Every spin has at least two neighbors of the same sign, since only loops and barbells are stable
under the effect of the dynamics.}\label{absorbing-state}
\end{center}
\end{figure}

The discrete time cellular automaton corresponding to the zero-temperature
case of Domany's stochastic Ising ferromagnet on the hexagonal lattice
can be considered as a simplified version of a continuous time Markov process 
where an independent (rate $1$) Poisson clock is assigned to each site
$x \in {\mathbb H}$, and the spin at site $x$ is updated
(with the same rule as in our discrete time dynamics)
when the corresponding clock rings. (In particular, they have the same stable
configurations -- see Figure~\ref{absorbing-state}.)
The percolation properties of the final state $\sigma^\infty$ of that
process  were studied, both rigorously and numerically, in~\cite{hn};
the results there (about critical exponents rather than the continuum
scaling limit)  strongly suggest that that dependent percolation model
is also in the same universality class as independent percolation.  
Similar stochastic processes on 
different types of lattices have been studied in various papers.
See, for example,~\cite{cdn, fss, gns, nns, ns1, ns2, ns3} 
for models on ${\mathbb Z}^d$
and~\cite{Howard} for a model on 
the homogeneous tree of degree three.
Such models are also discussed 
extensively in the physics literature, usually on ${\mathbb Z}^d$
(see, for example,~\cite{domany} and~\cite{lms}).
Numerical simulations have been done by Nienhius~\cite{Nienhuis} and
rigorous results for both the continuous and discrete dynamics
have been obtained in~\cite{cns}, including
a detailed analysis of the discrete time (synchronous) case.

Let us now describe in more detail the two deterministic cellular
automata (on ${\mathbb H}$ and ${\mathbb T}$).
Later, we will prove the equivalence of the percolation models
generated by those cellular automata. \\

\noindent {\bf Zero-temperature Domany model}

\noindent Consider the homogeneous ferromagnet
on the hexagonal lattice ${\mathbb H}$ with states denoted
by  $\sigma = \{ \sigma_x \}_{x \in {\mathbb H}}, \, \sigma_x = \pm 1$, 
and with (formal) Hamiltonian
\begin{equation} \label{Hamiltonian}
{\cal H}(\sigma) = - \sum_{ \langle x,y \rangle } \sigma_x \sigma_y ,
\end{equation}
where $\sum_{ \langle x,y \rangle }$ denotes the sum over all pairs of 
neighbor sites, each pair counted once.
We write ${\cal N}^{\mathbb H}(x)$ for the set of 
three neighbors of $x$, and indicate with
\begin{equation} \label{variation}
\Delta_x  {\cal H} (\sigma) =
2 \sum_{y \in {\cal N}^{\mathbb H}(x)} \sigma_x \sigma_y
\end{equation}
the change in the Hamiltonian when the spin $\sigma_x$ at site $x$ is 
flipped (i.e., changes sign).

\begin{figure}[!ht]
\begin{center}
\includegraphics[width=8cm]{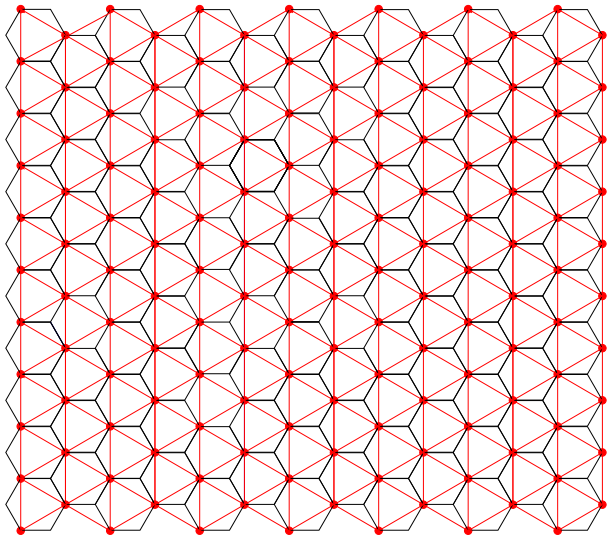}
\caption{A star-triangle transformation.}\label{star-tri}
\end{center}
\end{figure}

The hexagonal lattice ${\mathbb H}$ is partitioned
into two subsets ${\cal A}$ and ${\cal B}$
in such a way that all three neighbors of a site $x$ in
${\cal A}$ (resp., ${\cal B}$) are in ${\cal B}$ (resp., ${\cal A}$).
By joining two sites of ${\cal A}$ whenever they are next-nearest neighbors
in the hexagonal lattice (two steps away from each other),
we get a triangular lattice (the same with ${\cal B}$)
(see figure~\ref{star-tri}).
The synchronous dynamics is such that all the sites in the sublattice
${\cal A}$ (resp., ${\cal B}$) are updated simultaneously.

We now define the discrete time Markov process
$\sigma^n, \, n \in {\mathbb N}$, with 
state space ${\cal S}_{\mathbb H} = \{ -1,+1 \}^{\mathbb H}$,
which is the zero temperature 
limit of a model of Domany~\cite{domany}, as follows:
\begin{itemize}
\item The initial state $\sigma^0$ is chosen 
from a symmetric Bernoulli product measure.
\item At odd times $n = 1, 3, \dots$, 
the spins in the sublattice ${\cal A}$ are updated according to the
following rule: $\sigma_x, \, x \in {\cal A}$, 
is flipped if and only if $\Delta_x  {\cal H} (\sigma)<0$.
\item At even times $n = 2, 4, \dots$, 
the spins in the sublattice ${\cal B}$ are updated according to the
same rule as for those of the sublattice ${\cal A}$.
\end{itemize}

\noindent {\bf Cellular automaton on ${\mathbb T}$}

\noindent We define here a deterministic cellular automaton Q
on the triangular lattice ${\mathbb T}$, with random initial state
chosen by assigning value $+1$ or $-1$ independently, with equal probability,
to each site of ${\mathbb T}$.

Given some site $\bar{x} \in {\mathbb T}$, group its six ${\mathbb T}$-neighbors 
$y$ in three disjoint pairs $\{ y_1^{\bar{x}}, y_2^{\bar{x}} \}$,
$\{ y_3^{\bar{x}}, y_4^{\bar{x}} \}$, $\{ y_5^{\bar{x}}, y_6^{\bar{x}} \}$, so that
$y_1^{\bar{x}}$ and $y_2^{\bar{x}}$ are ${\mathbb T}$-neighbors,
and so on for the other two pairs.
Translate this construction to all sites $x \in {\mathbb T}$, thus producing
three pairs of sites  $\{ y_1^x, y_2^x \}$, $\{ y_3^x, y_4^x \}$,
$\{ y_5^x, y_6^x \}$ associated to each site $x \in {\mathbb T}$.
(Note that this construction does not need to specify how
${\mathbb T}$ is embedded in ${\mathbb R}^2$.)
Site $x$ is updated at times $m = 1, 2, \ldots$ according to the following rule:
the spin at site $x$ is 
changed from $\sigma_x$ to $- \sigma_x$ if and only if 
at least two of its pairs of
neighbors have the same sign and this sign is $- \sigma_x$.

\boldmath
\subsection{Equivalence between the models on ${\mathbb H}$ and ${\mathbb T}$}
\unboldmath

We show here how the models on the hexagonal and on the triangular lattice
are related through a \emph{star-triangle transformation}.
More precisely, we will show that
the dynamics on the triangular lattice ${\mathbb T}$ is equivalent to
the alternating sublattice dynamics on the hexagonal lattice ${\mathbb H}$
when restricted to the sublattice ${\cal B}$ for even times $n=2m$.

To see this, start with ${\mathbb T}$ and construct an
hexagonal lattice ${\mathbb H}'$ by means of a star-triangle transformation
(see, for example, p. 335 of \cite{grimmett})
such that a site is added at the center of each of the triangles
$(x, y_1^x, y_2^x), (x, y_3^x, y_4^x)$, and $(x, y_5^x, y_6^x)$
(the sites $y_i^x$ are defined in the previous subsection).
${\mathbb H}'$ may be partitioned into two triangular sublattices ${\cal A}'$
and ${\cal B}'$ with ${\cal B}' = {\mathbb T}$.
One can now see that the dynamics on ${\mathbb T}$ for $m=1,2,\ldots$
and the alternating sublattice dynamics on ${\mathbb H}'$ restricted to
${\cal B}'$ for even times $n=2m$ are the same.

An immediate consequence of this equivalence between the two cellular automata
is that the two families of percolation models that they produce are also
equivalent in an obvious way through a star-triangle transformation.
To be more precise, the percolation models defined on ${\mathbb T}$ by $Q$
for times $m=1,2,\ldots$ are the same as those defined on ${\cal B}$ by
the zero-temperature Domany model for even times $n=2m$.

\subsection{Results for the zero-temperature Domany model}

In this section we explain how the results obtained
for the percolation models $\mu^n_{\delta}$ generated by the cellular
automaton T are also valid for the percolation models defined by
the zero-temperature Domany model.
Because of the results of the previous section, we can actually
consider the percolation models $\tilde\mu^n_{\delta}$ on
$\delta{\mathbb T}$ generated by Q, for which we have results
analogous to Theorem~\ref{thm} and Corollary~\ref{cor}.
Results equivalent to Propositions~\ref{fixation},~\ref{perc}
and~\ref{mean_c_s} are contained in~\cite{cns}.
The proof of the main theorem (i.e., the analogue of Theorem~\ref{thm})
is basically the same as for the models
generated by T, so we just point out the differences.
We follow here the setup and notation of~\cite{cns}, but give all
the relevant definitions in order to make this paper self-consistent.

Let us consider a loop $\Gamma$ in the 
triangular sublattice ${\cal B}$, written as an ordered
sequence of sites $(y_0, y_1, \dots, y_k)$ 
with $k \geq 3$, which are distinct except that $y_k=y_0$.
For $i=1, \dots, k$, let $\zeta_i$ be the 
unique site in ${\cal A}$ that is an ${\mathbb H}$-neighbor
of both $y_{i-1}$ and $y_i$.
We call $\Gamma$ an {\bf s-loop} if $\zeta_1, \dots, \zeta_k$ are all distinct.
Similarly, a (site-self avoiding)
path $(y_0, y_1, \dots, y_k)$ in ${\cal B}$, between
$y_0$ and $y_k$, is called an {\bf s-path} if
$\zeta_1, \dots, \zeta_k$ are all distinct. 
Notice that any path in ${\cal B}$ between $y$ and $y'$
(seen as a collection of sites) contains an s-path between $y$ and $y'$.
An s-loop of constant sign is stable 
for the dynamics since at the next update of ${\cal A}$
the presence of the constant sign 
s-loop in ${\cal B}$ will produce a stable loop of that sign
in the hexagonal lattice.
Similarly an s-path of constant sign between $y$ and $y'$ will be stable 
if $y$ and $y'$ are stable --- e.g., if they each belong to an s-loop.
A triangular loop $x_1, x_2, x_3 \in {\cal B}$ with a common
${\mathbb H}$-neighbor $\zeta \in {\cal A}$ is called a {\bf star}; it is not 
an s-loop.
A triangular loop in ${\cal B}$ 
that is not a star is an s-loop and will be called an {\bf antistar},
while any loop in ${\cal B}$ 
that contains more than three sites contains an s-loop.

With these definitions, the proof of the main theorem for the zero-temperature
Domany model is the same as that presented in Section~\ref{proofs} for our
prototypical model, with the role of loops in the original proof played
here by s-loops (in particular antistars), and that of m-paths by s-paths
(see Section~\ref{proofs} below).

\boldmath
\subsection{An amusing further example: totally synchronous dynamics on ${\mathbb H}$}
\unboldmath

The model that we consider here corresponds to the zero-temperature Glauber
dynamics on the hexagonal lattice with \emph{all} sites updated simultaneously
at discrete times, \emph{without} alternating between two subsets (contrary to
the case of the zero-temperature Domany model), with initial configuration
$\sigma^0$ chosen according to a symmetric Bernoulli product measure.

Let us partition ${\mathbb H}$ in two subsets ${\cal A}$ and ${\cal B}$
as before and define the family of percolation models $\{ \bar\mu^m_{\cal A},
m = 0,1, \ldots \}$, where $\bar\mu^m_{\cal A}$ is the distribution of
$\sigma^{2m}$ (at even times) restricted to the subset ${\cal A}$
(naturally endowed with a triangular lattice structure, so that
$\bar\mu^0_{\cal A}$ is the distribution of critical site percolation).

The main difference consists in the fact that $\sigma^n$ does not fixate
as $n \to \infty$ since there is a positive density of spins that flip
infinitely many times.
To see this, consider a loop $\Gamma$ in ${\mathbb H}$ containing an even
number of sites and such that at time zero its spins are alternately plus
and minus.
At any time $n$, every spin in $\Gamma$ has two neighbors of opposite sign
and will therefore flip at the next update.
Thus, the spins in $\Gamma$ never stop flipping.

However, there is a simple observation that tremendously simplifies
the analysis of this totally synchronous dynamics.
Namely, that if we restrict attention to sublattice ${\cal A}$ at odd times
$n = 1, 3, 5, \ldots$ and sublattice ${\cal B}$ at even times $n = 0, 2, 4, \ldots$,
the dynamics is \emph{identical} to the zero-temperature Domany dynamics
discussed above (let us call this $\sigma^n_a$).
On the other hand, if we instead observe ${\cal A}$ at even times and ${\cal B}$
at odd times, this is identical to an alternative zero-temperature Domany type
dynamics $\sigma^n_b$ but with the first (and third and $\ldots$) update on
${\cal B}$ rather than ${\cal A}$.
Furthermore, $(\sigma^n_a)^{\infty}_{n=0}$ and $(\sigma^n_b)^{\infty}_{n=0}$
are completely independent of each other.
We conclude that there are two distinct limits $\sigma_a^{\infty}$ and
$\sigma_b^{\infty}$ (independent of each other) and that the scaling limit
of $\sigma^n$ restricted to either ${\cal A}$ or ${\cal B}$ is the same
as for independent critical percolation on ${\mathbb T}$.
But it appears that for any $n$, $\sigma^n$ on all of ${\mathbb H}$ should
be subcritical and thus have a trivial scaling limit.

\end{document}